\documentclass[11pt]{article}

\usepackage{amssymb}
\usepackage[mathscr]{eucal}
\usepackage[all]{xy}
\usepackage{cancel}
\usepackage{graphicx}

\usepackage{amsmath}

\begin{document}

\noindent {\bf On Yiu's Equilateral Triangles Associated with a Kiepert Hyperbola}\\

\author*{\centerline {Cherng-tiao Perng}}
\vskip 0.2in

\let\thefootnote\relax\footnotetext{In memory of my father, Hsin-Hsi Perng, and my colleague, Dr. Boyd Coan}

\noindent {\bf Abstract.} In 2014, Paul Yiu constructed two equilateral triangles inscribed in a Kiepert hyperbola associated with a reference triangle. It was asserted that each of the equilateral triangles is triply perspective to the reference triangle, and in each case, the corresponding three perspectors are collinear. In this note, we give a proof of his assertions. Furthermore as an analogue of Lemoine's problem, we formulated and answered the question about how to recover the reference triangle given a Kiepert hyperbola, one of the two Fermat points and one vertex of the reference triangle.\\

\section{Introduction}

\noindent Given a scalene triangle $ABC$, one has the first Fermat point $F_1$ and the second Fermat point $F_2$. The \emph{Kiepert hyperbola} $\mathcal{K}$ we need here is the unique conic passing through the five points $A,B,C,F_1$ and $F_2$. Since it is well known that $\mathcal{K}$ must also pass through the centroid $M$ of $ABC$, it is convenient to construct $\mathcal{K}$ from the five points $A,B,C,F_1$ and $M$. We will call $ABC$ the \emph{reference triangle} of $\mathcal{K}$. Yiu's problem is stated as follows.\\

\noindent {\bf Theorem 1.} Let $\mathcal{K}$ be a Kiepert hyperbola with reference triangle $ABC$. Let $\mathcal{C}$ (resp. $\mathcal{C'}$) be the circle centered at $F_2$ (resp. $F_1$) with radius $F_2F_1$. The circle $\mathcal{C}$ (resp. $\mathcal{C'}$) intersects $\mathcal{K}$ at three points $P,Q,R$ (resp. $P',Q',R'$) other than $F_1$ (resp. $F_2$). Then with suitable order (clockwise or counterclockwise), one has

\noindent (a) Triangle $PQR$ (resp. $P'Q'R'$) is equilateral (See Figure 1);\\
\noindent (b) Triangle $PQR$ (resp. $P'R'Q'$) is triply perspective to triangle $ABC$, in the sense that $PA,QB,RC$ meet at a common point $P_1$; $PB,QC,RA$ meet at a common point $P_2$; $PC,QA,RB$ meet at a common point $P_3$; and moreover\\
\noindent (c) The three perspectors $P_1,P_2,P_3$ in (b) are collinear (in a line which also passes through the other three perspectors constructed from the triangle $P'R'Q'$).\\

\begin{figure}[!htbp]
\centering
\includegraphics[width=0.5\textwidth]{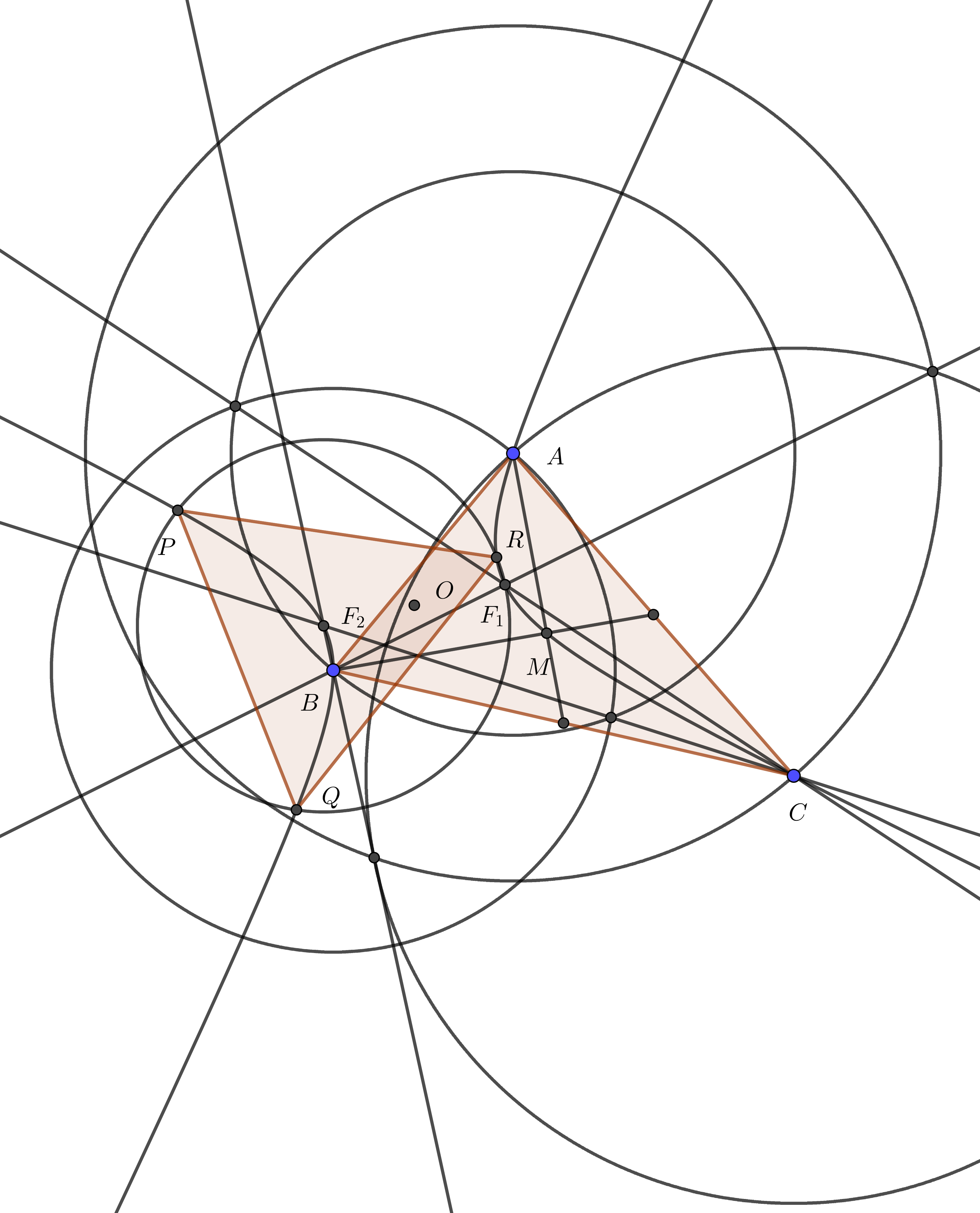}
\caption{An equilateral triangle inscribed in a Kiepert hyperbola}
\end{figure}

\noindent The author first learned Paul Yiu's problem in July 2018 from the Mathoverflow website ([13]; See also [2]), where Dao Thanh Oai posted and promoted Yiu's problem. According to Dao, this result is comparable to the case of Morley triangles or Napoleon triangles, and Dao placed a great importance to this construction as Kiepert hyperbola is a very special conic which passes through many triangle centers. The author has applied two methods to solve Yiu's problem. The first method (of which the details will not be included in this note) involves direct computation by introducing the slope variable $m=(y-y_2)/(x-x_2)$, where $(x_2,y_2)$ is the second Fermat point. By elimination, $m$ will need to satisfy a cubic equation $p(m)=0$. Then by manipulations with symmetric functions on the roots of $p(m)$, one can prove all the statements of the problem. The second method is much simpler and gives a construction of the reference triangle, assuming that the Kiepert hyperbola, one of the two Fermat points, and one of the vertices of the reference triangle are given. Before we go into the proofs or constructions in Section 3, we briefly recall some basic facts and definitions in the next section. In Section 4, we state and prove a result which is supposed to be known more than a hundred years ago.\\

\section{Preliminaries}
\noindent We collect some definitions and known facts here.\\

\noindent {\bf 2.1.} The first Fermat point $F_1$ (or the first isogonic center) of a triangle is constructed by erecting an equilateral triangle from the outside of each side of the triangle: the lines joining the outer vertex to the opposite vertex of the triangle concur at $F_1$.\\

\noindent {\bf 2.2.} Similarly, the second Fermat point $F_2$ (or the second isogonic center) of a triangle is constructed by erecting an equilateral triangle from the inside of each side of the triangle.\\

\noindent {\bf 2.3.} The Nine Point Circle (or the Feuerbach circle) for a triangle is the circle that passes through the following nine points: the midpoint of each side of the triangle, the foot of each altitude, and the midpoint of the line segment from each vertex of the triangle to the orthocenter.\\

\noindent {\bf 2.4.} The Kiepert hyperbola defined in the beginning introduction is a rectangular hyperbola (meaning that the asymptotes are perpendicular to each other). It is known that for a triangle inscribed in a rectangular hyperbola, the orthocenter of the triangle lies on the rectangular hyperbola.\\

\noindent {\bf 2.5.} A version of the Feuerbach Conic Theorem we will use is the following: Let $ABC$ be a triangle inscribed in a rectangular hyperbola. Then the Feuerbach circle of $ABC$ passes through the center of the rectangular hyperbola. Furthermore the center of the Kiepert hyperbola associated with $ABC$ is midway between the two Fermat points (isogonic centers) ([1],[3]).\\

\noindent Historically, Kiepert hyperbola was introduced by Ludwig Kipert to solve Lemoine's problem (Question 864 in Nouvelles Annales de Math\'{e}matiques, Series 2, Volume 7 (1868), p. 191). The question published by Lemoine asks about the following construction problem: Given one vertex of each of the equilateral triangles placed on the sides of a triangle, construct the original triangle (Wikipedia [12]). Kiepert solved the problem and since then, Kiepert hyperbola has been widely studied. As an analogue of Lemoine's problem and in the setting of Yiu's problem, we may pose the following question: Given a Kiepert hyperbola, one of the two Fermat points and one vertex of the reference triangle, construct the other two vertices of the reference triangle. The question is answered directly by the main results of our note.\\

\noindent {\bf 2.6.} Any two conics are related by a projective collineation: Any three distinct points on the first conic can be made to correspond to any three distinct points of the second (Ex.4 on page 79 of [9], or Theorem 6.4.1 of [10]).\\

\noindent {\bf Lemma 2.7.} Let $ABC$ and $A'B'C'$ be two triangles inscribed in a conic $\mathcal{K}$ such that $ABC$ is perspective with $A'B'C'$ with perspector $P$, and $ABC$ is perspective with $B'C'A'$ with perspector $Q$. Then there exists $R\in PQ$ such that $ABC$ is perspective with $C'A'B'$ with perspector $R$.\\

\noindent {\it Proof.} The existence of $R$ such that $ABC$ is perspective with $C'A'B'$ is guaranteed since it is well known that doubly perspective implies triply perspective. It suffices to show that the three perspectors $P,Q,R$ are collinear. Applying Pascal's theorem to the ordered lists $\{C',B',A'\}$ and $\{A,B,C\}$ shows that $Q,R$ and $P$ are collinear. $\Box$\\

\noindent {\bf Definition.} Let $ABC$ be a triangle inscribed in a conic $\mathcal{K}$. Let the line tangent to $\mathcal{K}$ at $A$ meet $BC$ at $P$, the line tangent to $\mathcal{K}$ at $B$ meet $CA$ at $Q$, and the line tangent to $\mathcal{K}$ at $C$ meet $AB$ at $R$. Then by a limiting case of Pascal's theorem, it can be shown that $P,Q,R$ are collinear. We call this line the \emph{Hessian line} of $ABC$ with respect to $\mathcal{K}$ (cf. p.130 of [11]).\\

\noindent {\bf Lemma 2.8.} Let $ABC$ and $A'B'C'$ be two triply perspective triangles inscribed in a conic. Then the line passing through the three perspectors (see Lemma 2.7) coincide with the Hessian line of $ABC$ with respect to the conic.\\

\noindent {\it Proof.} Let the tangent line at $A$ (resp. $B$) meet $BC$ (resp. $CA$) at $P_1$ (resp. $P_2$). Let $ABC$ be perspective with $A'B'C'$ (resp. $C'A'B'$) of perspector $Q_1$ (resp. $Q_3$). Applying Pascal theorem to the ordered lists $\{A,B,C'\}$ and $\{C,A,A'\}$ (resp. $\{B,C,A'\}$ and $\{A,B,B'\}$) shows that $P_1$ (resp. $P_2$) lies on $Q_1Q_3$. Hence lines $P_1P_2$ and $Q_1Q_3$ coincide. $\Box$\\

\section{Proof of Theorem 1}

\noindent {\it Proof of (a) in Theorem 1.} By construction, $P,Q,R,F_1$ lie on a circle of center $F_2$ and radius $F_2F_1$. Since $\mathcal{K}$ is a rectangular conic containing the points $P,Q,R,F_2$, we can apply Feuerbach's theorem (2.5) to the triangles $F_2PQ,F_2QR,F_2RP$ and $PQR$ which are inscribed in $\mathcal{K}$. Let $D,E,F$ be the midpoints of $PQ,QR,RP,$ respectively, and let $U,V,W$ be the midpoints of $F_2P,F_2Q,F_2R$, respectively. Note that by construction, $$F_2U=F_2V=F_2W=\frac 1 2F_2P=\frac 1 2F_2F_1=F_2O,$$ where $O$ is both the center of the Kiepert hyperbola and the midpoint of $F_2F_1$ (See (2.5)). By Feuerbach's theorem applied to the four triangles and the conic $\mathcal{K}$ with center $O$, it follows that $$D,U,V,O~{\rm are~ concyclic~ with~ center ~}F_2,$$
     $$E,V,W,O~{\rm are~ concyclic~ with~ center ~}F_2,$$
     $$F,W,U,O~{\rm are~ concyclic~ with~ center ~}F_2,$$ and
     $$D,E,F,O~{\rm are~ concyclic~ with~ center ~}F_2,$$ where the last statement follows from the first three. Now since $F_2$ is circumcenter of $PQR$, the points $D,E,F$ are perpendicular feet from $F_2$ to the sides of the triangle $PQR$. Pythagorean Theorem shows now that $PQR$ is equilateral. $\Box$\\

\noindent In order to prove the other parts of Theorem 1, we reverse the process by looking into the following situation. Without loss of generality, the rest assertions of Theorem 1 follow from the following theorem.\\

\noindent {\bf Theorem 2.} Let $F_2=(-1,0),F_1=(1,0)$ and let $S_1$ (resp. $S_2$) be a circle centered at $F_2$ (resp. $F_1$) with radius $F_2F_1=2$. Let $PQR$ be an equilateral triangle inscribed in $S_1$ and let $P'Q'R'$ (an equilateral triangle inscribed in $S_2$) be the reflection of $PQR$ with respect to the center $(0,0)$ of $F_2F_1$ (See Figure 2, where triangle $P'Q'R'$ is not drawn for simplicity). Then\\

\noindent (a) The two triangles $PQR$ and $P'Q'R'$ are triply perspective with perspectors consisting of $(0,0)$ and the intersection of the two circles $S_1$ and $S_2$. Accordingly the perspectors lie on the radical axis $L$ of $S_1$ and $S_2$.\\

\noindent (b) There is a unique conic $\mathcal{K}$ passing through $P,Q,R,P',Q',R',F_2$ and $F_1$.\\

\noindent (c) For any point $V\in L$, let $P''$ be the other point of intersection of the line $PV$ and $\mathcal{K}$. Define $Q''$ and $R''$ in a similar way. Then $PQR$ is triply perspective with $P''Q''R''$, and the other two perspectors lie on $L$. Similarly, $P'Q'R'$ is triply perspective with $P''R''Q''$ (i.e. one needs to reverse the orientation).\\

\noindent (d) Let $P''Q''R''$ be the triangle constructed in (c). Then $\mathcal{K}$ is the Kiepert hyperbola through the two Fermat points $F_1$ and $F_2$ with respect to the reference triangle $P''Q''R''$.\\

\noindent (e) For two distinct triangles $P_1Q_1R_1$ and $P_2Q_2R_2$ constructed as in (c), $P_1Q_1R_1$ is triply perspective with $P_2R_2Q_2$. These are all triply perspective with $P'R'Q'$.\\

\noindent {\it Proof.} (a) It is easy to see that $S_1$ and $S_2$ intersect at $\{(0,2),(0,-2)\}$ and $L$ is the line with equation $x=0$. The result in (a) is easy to see either by synthetic method or by coordinates.\\

\noindent (b) Denoting $P$ (using rational parametrization of a circle) by $$P=\left(-1+\frac{2(1-t^2)}{1+t^2},\frac {4t}{1+t^2}\right)$$ and turning counterclockwise $120^\circ$ successively, we have $$Q=\left(\frac{-2(\sqrt{3}t+1)}{t^2+1},\frac{-(\sqrt{3}t^2+2t-\sqrt{3})}{t^2+1}\right)$$ and $$R=\left(\frac{2(\sqrt{3}t-1)}{t^2+1},\frac{\sqrt{3}t^2-2t-\sqrt{3}}{t^2+1}\right).$$ Then the unique conic $\mathcal{K}$ through $P,Q,R,F_1,F_2$ has the equation $$(3t^2-1)x^2+(2t^3-6t)xy+(1-3t^2)y^2+(1-3t^2)=0$$ and it is evident that $\mathcal{K}$ also passes through $P',Q',R'$, which are the reflections of $P,Q,R$ with respect to the origin.\\

\noindent (c) Letting $V=(0,y_0)$ be a point on $L$, we solve the alternative points on the conic and get that $$P''=\left(\frac{(3t^2-1)(y_0^2+1)}{(t^2+1)(2ty_0-y_0^2+1)},\frac{2(t^2y_0+2t-y_0)(ty_0-1)}{(t^2+1)(2ty_0-y_0^2+1)}\right),$$
$$Q''=\left(\frac{-2(3\sqrt{3}~t^3+3t^2-\sqrt{3}~t-1)(y_0^2+1)}{3t^2y_0^2+2\sqrt{3}~t^2y_0-3t^2+8ty_0-y_0^2+2\sqrt{3}~y_0+1)(t^2+1)},\right.$$
$$ -(\sqrt{3}~t^4y_0^2+6t^4y_0-2t^3y_0^2+3\sqrt{3}~t^4-8\sqrt{3}~t^2y_0^2+6t^3+4t^2y_0-10ty_0^2-4\sqrt{3}~t^2-\sqrt{3}~y_0^2$$
$$\left.-2t-2y_0+\sqrt{3})/((3t^2y_0^2+2\sqrt{3}~t^2y_0-3t^2+8ty_0-y_0^2+2\sqrt{3}~y_0+1)(t^2+1))\right)$$ and
$$R''=\left(\frac{2(3\sqrt{3}~t^3-3t^2-\sqrt{3}~t+1)(y_0^2+1)}{3t^2y_0^2-2\sqrt{3}~t^2y_0-3t^2+8ty_0-y_0^2-2\sqrt{3}~y_0+1)(t^2+1)},\right.$$
$$ (\sqrt{3}~t^4y_0^2-6t^4y_0+2t^3y_0^2+3\sqrt{3}~t^4-8\sqrt{3}~t^2y_0^2-6t^3-4t^2y_0+10ty_0^2-4\sqrt{3}~t^2-\sqrt{3}~y_0^2$$
$$\left.+2t+2y_0+\sqrt{3})/((3t^2y_0^2-2\sqrt{3}~t^2y_0-3t^2+8ty_0-y_0^2-2\sqrt{3}~y_0+1)(t^2+1))\right).$$ Now it is straightforward (by a computer algebra system) to check that $$PP'',QQ'',RR''~{\rm are~concurrent~at~}(0,y_0)~{\rm (by~construction)},$$
$$PQ'',QR'',RP''~{\rm are~concurrent~at~}(0,-(\sqrt{3}~y_0+3)/(3y_0-\sqrt{3})),$$ and
$$PR'',QP'',RQ''~{\rm are~concurrent~at~}(0,(\sqrt{3}~y_0-3)/(3y_0+\sqrt{3})).$$

Parts (d) and (e) are easily checked by a computer algebra system. $\Box$ \\

\noindent With Theorem 2, the proof of Theorem 1 is completed. $\Box$\\

\noindent In the spirit of Lemoine's problem, we may ask about how to recover the reference triangle used to construct the Kiepert hyperbola. We have the following result. \\

\begin{figure}
\centering
\includegraphics[width=0.5\textwidth]{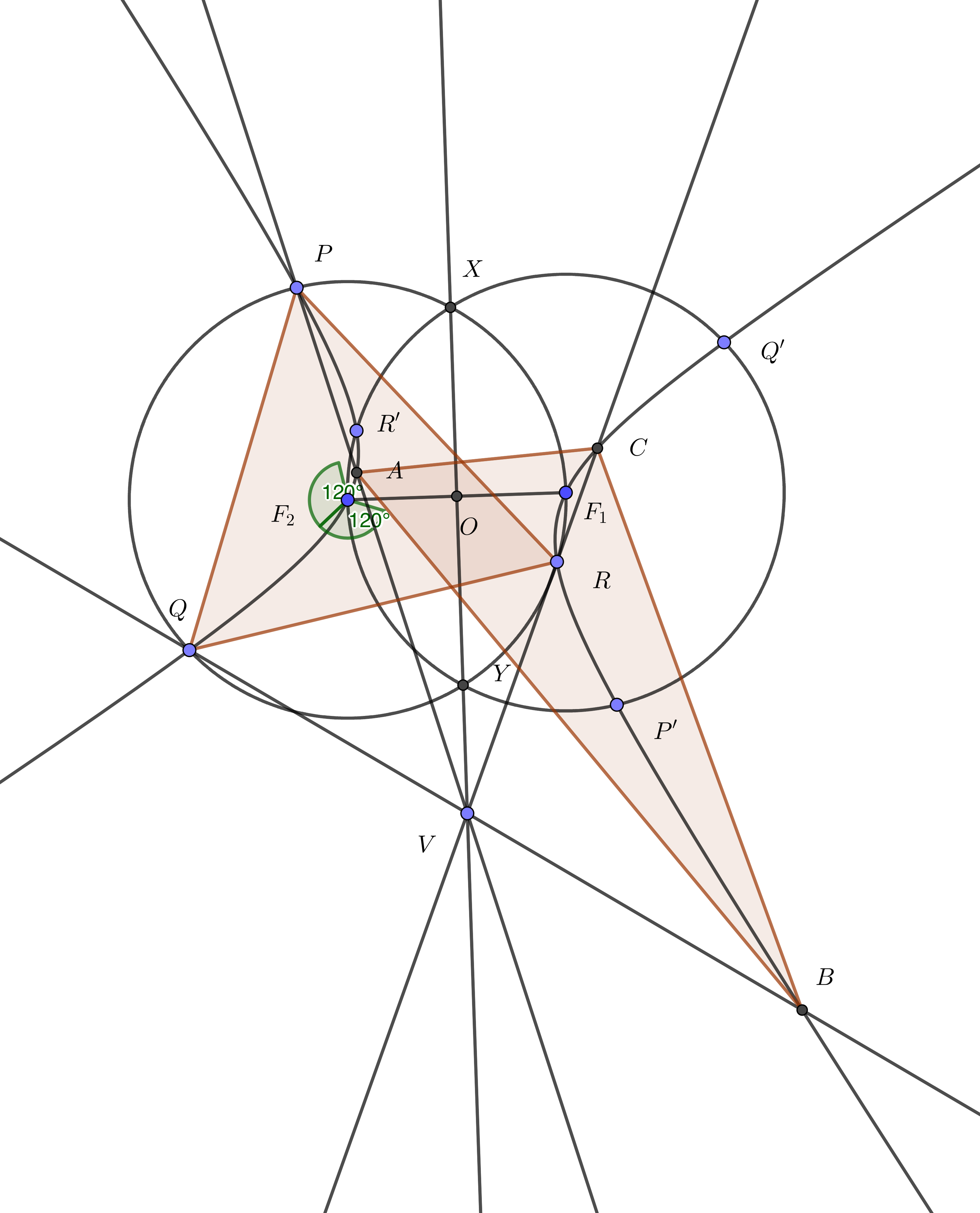}
\caption{Construction of a reference triangle}
\end{figure}

\noindent {\bf Corollary.} Given a Kiepert hyperbola $\mathcal{K}$, one of the two Fermat points and a vertex of the reference triangle on $\mathcal{K}$, we can recover the other two vertices on the reference triangle.\\

\noindent {\it Proof.} Note that the center $O$ of $\mathcal{K}$ can be constructed by using 2.5. Reflecting one Fermat point with respect to $O$ gives another Fermat point. With $\mathcal{K}$ and the two Fermat points, we follow the construction as in the statement of Theorem 2 and the proof for (a) of Theorem 1 to get an equilateral triangle $PQR$. Using the notations and the construction as in the Theorem 2, join any vertex of $PQR$, say $P$, with a vertex of the reference triangle, say $A$; line $PA$ intersects $L$ at a point $V$ (See Figure 2). Then follow the steps in (c) to recover $B$ and $C$. $\Box$\\

\section{A Generalization}

\noindent The following theorem gives a generalization of some results in Theorem 2 without computation by drawing facts from projective geometry. Similar or more general results (with different proofs) were known (cf. [4],[5],[6],[7],[8]).\\

\noindent {\bf Theorem 3.} Let $ABC$ and $A'B'C'$ be two triangles inscribed in a conic $\mathcal{K}$ such that $ABC$ is perspective with $A'B'C'$ with perspector $P$, and $ABC$ is perspective with $B'C'A'$ with perspector $Q$. Then\\

\noindent (a) There exists $R\in PQ$ such that $ABC$ is perspective with $C'A'B'$ with perspector $R$, and $P,Q,R$ are collinear.\\

\noindent (b) If $A''B''C''$ is another triangle inscribed in $\mathcal{K}$ and $ABC$ is perspective with $A''B''C''$ with perspector $S$ for some $S\in PQ$, then there exist $T,U\in PQ$ such that $ABC$ is perspective with $B''C''A''$ (resp. $C''A''B''$) with perspector $T$ (resp. $U$).\\

\noindent {\it Proof.} Part (a) is already proven in Lemma 2.7. To prove (b), we apply a projective collineation to $ABC$ on $\mathcal{K}$ and send it to the vertices of an equilateral triangle on a circle. For simplicity we keep using the same letters (so assume $\mathcal{K}$ is a circle now). Clearly the Hessian line of $ABC$ with respect to the circle is the line at infinity. By Lemma 2.8, $P,Q,R$ lie on the line at infinity. Since $ABC$ and $A'B'C'$ are triply perspective with perspectors on the line at infinity, it is clear that $A'B'C'$ is an equilateral triangle. Note that after the collineation, $A''B''C''$ is a triangle on the circle which is perspective with $ABC$ with a perspector $S$ on the line at infinity, i.e. in the affine plane $AA'', BB'',CC''$ are parallel. Necessarily, $A''B''C''$ is an equilateral triangle (in reverse orientation as $ABC$). But now it is easy to see that $ABC$ and $A''B''C''$ are triply perspective, again with perspectors lying on the line at infinity. Mapping back by the inverse projective collineation, we have proven the required result. $\Box$\\

\noindent {\bf Acknowledgment.} We use the computer algebra system SAGE [14] in our computation.\\

\noindent {\bf References.}\\

\noindent [1] H. Martyn Cundy, Feuerbach's Theorem and the Rectangular Hyperbola, The Mathematical Gazette, Vol. 43, No. 343 (Feb, 1959), p.21-23.\\

\noindent [2] Dao Thanh Oai, Some New Equilateral Triangles in a Plane Geometry, Global Journal of Advanced Research on Classical and Modern Geometries, Vol. 7 (2018), Issue 2, p.73-91.\\

\noindent [3] R.H. Eddy and R. Fritsch,  The Conics of Ludwig Kiepert: A Comprehensive Lesson in the Geometry of the Triangle, Mathematics Magazine, Vol. 67, No. 3, June 1994, p.188-205.\\

\noindent [4] L. Ripert, Etude sur des groupes de triangles trihomologiques inscrits ou circonscrits \`{a} une m\`{e}me conique ou \`{a} des familles de coniques, AFAS, Tom 2, 1900, pp. 112-133.\\

\noindent [5] L. Ripert, Sur les triangles trihomologiques, Mathesis, 1900, p.226.\\

\noindent [6] J.A. Third, Triangles triply in perspective, Proceedings of the Edinburgh Mathematical Society, Vol. 19, Feb. 1900, pp.10-22, Published online by Cambridge University Press on 20 January 2009.\\

\noindent [7] Charles M'Leord and William P. Milne, Triangles Triply in Perspective, Vol. 28, Feb. 1909, pp.148-151, Published online by Cambridge University Press on 20 January 2009.\\

\noindent [8] F.G. Stockton, A set of triply perspective triangles associated with projective triad, AMM Vol. 62, No. 7, (1955), pp.41-51.\\

\noindent [9] H.S.M. Coxeter, Projective Geometry, First Edition, 1964, Blaisdell Publishing Company.\\

\noindent [10] G. Glaeser, H. Stachel, and B. Odehnal, The Universe of Conics: From the ancient Greeks to 21st century developments, Springer Spektrum, 2016.\\

\noindent [11] D.J. Struik, Lectures on Analytic and Projective Projective Geometry, Addison-Wesley Publishing Company, Inc., 1953.\\

\noindent [12] URL: https://en.wikipedia.org/wiki/Lemoine$\%27s\_$problem (last retrieved on Sep 2, 2019).\\

\noindent [13] URL: https://mathoverflow.net/questions/303722/yius-equilateral-triangle-triplet-points (last retrieved on Sep 2, 2019).\\

\noindent [14] SageMath, the Sage Mathematics Software System (Version 6.10), The Sage Developers, https://www.sagemath.org.\\

\end{document}